\documentclass[11pt]{amsart}

\usepackage{accents}
\usepackage{amsmath}
\usepackage{amsthm}
\usepackage{amsfonts}
\usepackage{bbm}
\usepackage{amssymb}
\usepackage{mathbbol}
\usepackage{tikz}
\usepackage{mathrsfs}

\usepackage{graphicx}
\usepackage{txfonts}
  \usepackage[figuresright]{rotating}
  \usepackage{floatpag} 
 \usepackage[all]{xy} 
 
 \usetikzlibrary{matrix,arrows,decorations.markings}
\newcommand{\Q}{\mathbb{Q}}
\newcommand{\Z}{\mathbb{Z}}
\newcommand{\C}{\mathbb{C}}

\newcommand{\mf[1]}{\mathfrak{#1}}

\newcommand{\g}{\mathfrak{g}}

\newcommand{\h}{\mathfrak{h}}

\DeclareMathOperator{\supp}{supp}
\DeclareMathOperator{\nil}{nil}
\DeclareMathOperator{\rank}{rank}

\DeclareMathOperator{\ad}{ad}

\newtheorem{Lem}{Lemma}[section]
\newtheorem{Thm}[Lem]{Theorem}
\newtheorem{Def}[Lem]{Definition}

\newtheorem{Rmk}[Lem]{Remark}

\title[ad-locally nilpotent elements of GKM]{Subalgebra generated by ad-locally nilpotent elements of Borcherds Generalized Kac-Moody Lie algebras} 
\author{Shrawan Kumar}

\address{Department of Mathematics, University of North Carolina, Chapel Hill, NC 27599-3250, USA}
\email{shrawan@email.unc.edu} 
\begin{document}

\maketitle

\section{Abstract} We determine the Lie subalgebra $\g_{\nil}$ of a Borcherds symmetrizable generalized Kac-Moody Lie 
 algebra $\g$ generated by $\ad$-locally nilpotent elements and show that it is `essentially' the same as the Levi subalgebra of $\g$ with its simple roots precisely the real simple roots of $\g$. 
\section{Introduction}
Let $\g= \g(A)$ be the  symmetrizable Generalized Kac-Moody
(GKM)  algebra associated to a $\ell\times \ell$ matrix $A$ (cf. Section 3). Let
$$\g^o_{\nil} := \{x\in \g: \,\,\text{$\ad x$ acts locally nilpotently on $\g$}\},$$
and let $\g_{\nil} \subset \g$ be the Lie subalgebra generated by $\g^o_{\nil} $.  Then, we prove the following theorem
(cf. Theorem \ref{thm1}):

\vskip1ex
\noindent
{\bf Theorem.}
 Let  $\g= \mf[g](A)$ be as above, where $\ell \geq 2$ and $A$ is indecomposable, i.e., the corresponding Dynkin diagram is connected. Then,
$$\g'(B) \subset \g_{\nil}\subset \g'(B) +\h,$$
where $B\subset A$ is the submatrix parameterized by those $i$ such that $a_{i,i}=2$
and $\g'(B)$ is the derived subalgebra of $\g(B)$.
\vskip1ex
As shown in Remark \ref{remark2}, the assumption $\ell \geq 2$ in the above theorem  is necessary in general. 

\vskip2ex

\noindent
{\bf Acknowledgements:} We thank V. Kac for providing some references on   generalized Kac-Moody Lie 
 algebras. We also acknowledge partial support from the NSF grant number DMS-1802328.
\section{Basic Definition}

In this section, we recall the definition of Borcherds Generalized Kac-Moody Lie algebras $\mf[g]$ (for short GKM algebras).  For a more extensive treatment of $\mf[g]$ and its properties, see Chapters 1, 11 of \cite{Kac} and the papers \cite{Bo1} and \cite{Bo2}. 

\begin{Def} {\rm Let $A = (a_{i,j})$ be a $\ell\times \ell$ matrix (for $\ell \geq 1$) with real entries, satisfying the following properties;
\vskip1ex
(P1) either $a_{i,i}=2$ or  $a_{i,i}\leq 0$,
\vskip1ex
(P2) $a_{i,j}\leq 0$ if $i\neq j$, and $a_{i,j}\in \Z$ if $a_{i,i} =2$,
\vskip1ex
(P3) $a_{i,j}=0$ if and only if $a_{j,i}=0$.

\vskip1ex
Fix a {\it realization}  of $A$, which is a triple $(\h, \Pi, \Pi^\vee)$ consisting of a complex vector space $\h$, $\Pi= 
\{\alpha_1, \dots, \alpha_\ell \}\subset \h^*$ and $\Pi^\vee = \{\alpha_1^\vee, \dots, \alpha_\ell^\vee\}\subset \h$ are indexed subsets, satisfying the following three conditions:

\vskip1ex
(Q1) both sets $\Pi$ and $\Pi^\vee$ are linearly independent,
\vskip1ex
(Q2)\,$\alpha_j (\alpha_i^\vee) = a_{i,j}$, for  all $i,j$,
\vskip1ex
(Q3) $\ell- \rank A = \dim \h- \ell$.
\vskip1ex

By \cite{Kac}, Proposition 1.1, such a realization is unique up to an isomorphism of the triple. 

Now, the Borcherds Generalized Kac-Moody Lie algebra (for short GKM algebra) $\mf[g](A)$ is defined as the Lie algebra generated by 
$\{e_i, f_i, \h\}_{1\leq i\leq \ell}$ subject to the following relations:

\vskip1ex
(R1) $[e_i, f_j] = \delta{ij}\alpha^\vee_i$, for all $i$,
\vskip1ex
(R2) $[h, h']= 0$,  for all $h,h'\in \h$,
\vskip1ex
(R3) $[h, e_i]= \alpha_i(h)e_i; \,\, [h, f_i]= -\alpha_i(h)f_i$, for all $1\leq i\leq \ell$ and $h\in \h$,
\vskip1ex
(R4) $(\ad e_i)^{1-a_{i,j}}e_j=  (\ad f_i)^{1-a_{i,j}}f_j =0,\,\,$ if $a_{i,i}=2$ and $i\neq j$,
\vskip1ex
(R5) $[e_i, e_j] = [f_i, f_j]=0,\,\,$ if $a_{i,j}=0$.
\vskip1ex

 The matrix $A$ (or the Lie algebra $\g(A)$) is called {\it symmetrizable} if there exists an invertible diagonal matrix
 $D = \text{diag} (\epsilon_1, \dots, \epsilon_\ell)$ such that the matrix $DA$ is symmetric.}

\end{Def}
\section{Main theorem and its proof}
\begin{Thm} \label{thm1} Let  $\g= \mf[g](A)$ be the symmetrizable GKM algebra associated to a $\ell\times \ell$ matrix $A$ as in the last section. Assume further that $\ell \geq 2$ and $A$ is indecomposable, i.e., the corresponding Dynkin diagram is connected. Let
$$\g^o_{\nil} := \{x\in \g: \,\,\text{$\ad x$ acts locally nilpotently on $\g$}\},$$
and let $\g_{\nil} \subset \g$ be the Lie subalgebra generated by $\g^o_{\nil} $. Then,
$$\g'(B) \subset \g_{\nil}\subset \g'(B) +\h,$$
where $B\subset A$ is the submatrix parameterized by those $i$ such that $a_{ii}=2$, i.e., $\alpha_i$ is a real root
and $\g'(B)$ is the derived subalgebra of $\g(B)$.
\end{Thm}
\begin{proof} Consider the  $\Z$-gradation of $\g$ induced from a homomorphism $\theta: Q:=\oplus_i\,\Z\alpha_i\to \Z$. Then, for any $x\in \g^o_{\nil}, x_+(\theta) \in \g^o_{\nil}$, where $x_+(\theta)$ is the top degree component of $x$ in the $\Z$-gradation of $\g$ induced by $\theta$. To prove this, observe that for any $y\in \g_\alpha$ (where $\g_\alpha$ is the root space corresponding to the root $\alpha$ or $0$),
$$(\ad x)^n(y) = (\ad x_+(\theta))^n(y) + \,\text{lower degree terms}.$$
Similarly, for $x\in \g^o_{\nil}, x_-(\theta) \in \g^o_{\nil}$, where $x_-(\theta)$ is the lowest degree component of $x$.

Further, given any nonzero $x\in \g$, we can get a gradation $\theta_x: Q\to \Z$ as above (depending upon $x$)  such that all the homogeneous degree components of $x$ (under $\theta_x$)  belong to root spaces $\g_\beta$. To prove this, write $x= \sum_j \,x_{\beta_j}$, where $\beta_j$ are distinct roots or zero, $x_{\beta_j}\in \g_{\beta_j}$  and each $x_{\beta_j}\neq 0$. Consider the finite collection of weights: $\{\beta_j -\beta_k\}_{j\neq k} \subset \h^*$. Now, we can find a vector $\gamma \in \mathbb{Q}^\ell = Q\otimes_{\Z}\Q$ such that the standard dot product $(\cdot , \cdot)$ in $\mathbb{Q}^\ell$:
\begin{equation} \label{eqn1} \theta_\gamma(\beta_j-\beta_k):= (\beta_j-\beta_k, \gamma) \neq 0, \,\,\text{for any $j\neq k$}.
\end{equation}
To prove the above equation, consider the $\ell -1$-dimensional subspace $V_{j, k}\subset \mathbb{Q}^\ell$ (for any $j\neq k$) perpendicular to $\beta_j -\beta_k$. Since the collection $\{\beta_j -\beta_k\}_{j\neq k}$ is finite, we can find a vector $\gamma$ such that the equation \eqref{eqn1} is satisfied. We can further take $\gamma \in Q\simeq \Z^\ell$ by clearing the denominators. 

So, if $x\in \g^o_{\nil}$, then either $x$ belongs to the center $Z(\g)$ of $\g$ (cf. \cite{Kac}, Proposition 1.6) or the root component $x_\beta\in \g^o_{\nil}$ for some root $\beta$ ($\beta \neq 0$). Moreover, if some nonzero root component of $x$ belongs to the root space $\g_{\delta}$ such that $\delta$ contains an imaginary  simple root $\delta_p$ (i.e., with $a_{p,p}\leq 0$) with nonzero coefficient, we can assume that $x_\delta\in  \g^o_{\nil}$ (possibly with a different nonzero root component of $x$ corresponding to a root containing an imaginary simple root with nonzero coefficient). This is achieved by taking $\gamma$ as above but requiring $\theta_\gamma(\alpha_p)$ to be much larger for all the imaginary simple roots $\alpha_p$ as compared to the values  $\theta_\gamma(\alpha_q)$ for all the real simple roots $\alpha_q$ (i.e., those with $a_{q,q}=2$).

By using the Cartan involution $\omega$ of $\g$ (i.e., $\omega (e_i) = -f_i, \omega(f_i) = -e_i, \omega (h) =-h\, \forall h\in \h$),  if needed, we can further assume that  $\delta$ is a positive root. Write
$$\delta = \sum_{p}  (m_p\alpha_p) +  \sum_{q}  (n_q\alpha_q),\,\,\text{for $m_p, n_q\geq 0$},$$
where $\alpha_p$ (resp. $\alpha_q$) run over all the  imaginary (resp. real) simple roots. In particular, some $m_p>0$. By \cite{Kac}, Exercise 11.21, the support 
$\supp (\delta)$ is connected. Assume first that $\delta$ is not an imaginary  simple root. Further, taking some $W$-translate (where $W$ is the Weyl group of $\g$, cf. \cite{Kac}, $\S$11.13), we can assume that $\delta(\alpha_q^\vee) \leq 0$ for all the real simple coroots $\alpha_q^\vee$ (cf. \cite{Kac}, Identity 11.13.3). Now, with respect to the $W$-invariant symmetric bilinear form $\langle \cdot, \cdot \rangle$ on $\h^*$ (cf. \cite{Kac}, $\S$2.1),
\begin{align} \label{eqn0} \langle \delta, \delta\rangle &= \sum_p\,m_p\langle \delta, \alpha_p\rangle + \sum_q\, n_q\langle \delta, \alpha_q\rangle\notag\\
 &=\sum_q\, n_q\langle \delta, \alpha_q\rangle + \sum_{p, q}\,m_pn_q\langle \alpha_q, \alpha_p\rangle + \sum_{p, p'}\,
 m_p m_{p'} \langle \alpha_{p'}, \alpha_{p}\rangle,
 \end{align}
 where $\alpha_{p'}$ also runs over imaginary simple roots. 
 Now, by assumption, 
 \begin{equation} \label{eqn2} \langle \delta, \alpha_q\rangle \leq 0,\,\,\text{ for all the real simple roots}.
 \end{equation}  
 For any imaginary simple root $\alpha_p$ 
 and any  real simple root  $\alpha_q$, 
 we have
  \begin{equation} \label{eqn3} \langle \alpha_q, \alpha_p\rangle \leq 0,\,\,\text{since $a_{p,q}\leq 0$}. 
  \end{equation}
  Further, for imaginary simple roots $\alpha_p, \alpha_{p'}$, 
   \begin{equation} \label{eqn4} \langle \alpha_{p'}, \alpha_{p}\rangle \leq 0,\,\,\text{by \cite{Kac}, Identity 2.1.6}. 
  \end{equation}
  Observe that we can take the normalizing factor $\epsilon_i >0$ for each $1\leq i\leq \ell$ as can be seen from 
  the identity: 
  $$\epsilon_ia_{i,j} = \epsilon_j a_{j,i},\,\,\text{for all $1\leq i, j \leq \ell$},$$
   where the diagonal matrix $D = \text{diag}(\epsilon_1, \cdots, \epsilon_\ell)$ is such that DA is a symmetric matrix. Moreover, since there exists $p$ with $m_p\neq 0$ and since $\supp \delta$ is connected and $\delta$ is not a simple root, by \cite{Kac}, Identity 2.1.6, 
   \begin{align} \label{eqn5} &\langle \alpha_{p'}, \alpha_{p}\rangle <0,\,\,\text{for some $p'\neq p$ with $m_{p'}\neq 0$ and $\alpha_{p'}$ an imaginary simple root}\notag\\
   &\text{or $\langle\alpha_q, \alpha_p\rangle <0$ for some $q$ with $n_q\neq 0$ and $\alpha_q$ a real simple root}. 
   \end{align}
   Thus, combining the equations \eqref{eqn0} - \eqref{eqn5}, we get:
   $$\langle \delta, \delta\rangle <0.$$
   By \cite{Kac}, Corollary 9.12, $\oplus_{k>0} \,\g_{k\delta}$ is a free Lie algebra on a basis of the form $\oplus_{k>0} \,\g_{k\delta}^o$, where
   $$\g_{k\delta}^o:=\{x\in   \g_{k\delta}:\langle x, y\rangle =0 \forall y\text{ in the Lie subalgebra generated by}  \,\g_{-\delta},\g_{-2\delta}, \dots, \g_{-(k-1)\delta}\}.$$
             Observe next that $\g_{k\delta}\neq 0$ for any $k>0$ by \cite{Kac}, Identity 11.13.3. If $\g_{\delta}$ is one dimensional, then so is $\g_{-\delta}$ and hence $\g_{2\delta}^o\neq 0$. Thus, $\oplus_{k>0} \,\g_{k\delta}$ is a free Lie algebra on at least $2$ generators. If $\dim \g_\delta \geq 2$, then  $\oplus_{k>0} \,\g_{k\delta}$ is again a free Lie algebra on at least two generators    (since $\g_\delta^o=\g_\delta$).           
              Thus, $\ad(x_\delta)$ can not act locally nilpotently on $\oplus_{k>0} \,\g_{k\delta}$ 
       and hence on $\g$  (since the enveloping algebra of a free Lie algebra is the tensor algebra on the same generators and now use \cite{Ku}, Identity (3) of Definition 1.3.2).
       
       Now, if $\delta = \alpha_p$ is an imaginary simple root. Then, again $x_\delta = e_p$ can not act nilpotently on any $e_i$, $i\neq p$  such that $a_{i,p}\neq 0$. (This is where we have used the assumption that $A$ is indecomposable and $\ell \geq 2$.) To prove this, use \cite{Kac}, Identity 11.13.3 by observing that $(n\alpha_p+ \alpha_i)\in K$ for all $n\geq 2$ in the notation of loc cit. Thus, we conclude that any $x\in \g^o_{\nil}$  must be of the form $x\in \g(B)+\h$. Hence,
       $$ \g_{\nil} \subset \g'(B)+\h.$$
        Further, by \cite{Ku}, Lemma 1.3.3(a) and the defining relations of $\g(A)$, $e_i, f_i\in  \g^o_{\nil}$ for any real simple root $\alpha_i$. Thus, 
        $$ \g'(B) \subset \g_{\nil}.$$
       This proves the theorem
       
       \end{proof}

\begin{Rmk} \label{remark2} {\rm (a) Define 
$${\g'}^o_{\nil} := \{x\in \g': \,\,\text{$\ad x$ acts locally nilpotently on $\g'$}\}$$
and let ${\g'}_{\nil}\subset \g'$ be the Lie subalgebra generated by ${\g'}^o_{\nil} $. Then, by the same proof as above,
$$\g'(B) \subset \g'_{\nil}\subset (\g'(B) +\h)\cap \g'.$$

\vskip1ex

(b) It is easy to see that the above theorem remains true in the case $A$ is parameterized by $\mathbb{N}\times \mathbb{N}$. 
\vskip1ex

(c) For the $1\times 1$-matrix $A =(0)$, following \cite{Kac}, $\S$2.9, $\g(A)= \h \oplus \C e_1 \oplus \C f_1$, where $\h =\C\alpha_1^\vee \oplus \C d$  and $[e_1, f_1] =\alpha_1^\vee, [\alpha_1^\vee, \g]=0, [d, e_1] =e_1, [d, f_1]= -f_1.$ Thus,  in this case, $ \g_{\nil}= \g'$.
Hence, the assumption $\ell \geq 2$ in the above theorem is necessary in general.
\vskip1ex

(d) One interesting consequence of the above theorem is that the only connected `reasonable' group attached to a GKM algebra $\g(A)$ is the one coming from its subalgebra $\g(B)$ (up to an $\h$-factor).}
\end{Rmk}

\end{document}